\documentclass{IEEEtaes}

\usepackage{color,array,amsthm}
\usepackage{graphicx}
\usepackage{amsmath}
\usepackage{tikz}
\usepackage{physics}
\usepackage{amssymb}
\usepackage{bm}
\usepackage{algorithm}
\usepackage{algpseudocode}
\usepackage{hyperref}
\usepackage{mathtools}
\usepackage{tabularx}
\usepackage{relsize}
\usepackage{subcaption} 

\usepackage[numbers]{natbib}
\usepackage{thmtools}
\usepackage{thm-restate}
\usepackage{cleveref}

\usepackage{multicol}

\usepackage{amsthm}

\jvol{XX}
\jnum{XX}
\jmonth{XXXXX}
\paper{1234567}
\pubyear{2024}
\doiinfo{TAES.2024.Doi Number}

\begin{document}

\title{Generalizing Trilateration: Approximate Maximum Likelihood Estimator for Initial Orbit Determination in Low-Earth Orbit}

\author{RICARDO FERREIRA}
\affil{NOVA School of Science and Technology, Caparica, Portugal} 

\author{FILIPA VALDEIRA}
\affil{NOVA School of Science and Technology, Caparica, Portugal}  

\author{MARTA GUIMAR\~{A}ES}
\affil{Neuraspace, Coimbra, Portugal}

\author{CL\'{A}UDIA SOARES}
\affil{NOVA School of Science and Technology, Caparica, Portugal} 


\receiveddate{Manuscript received XXXXX 00, 0000; revised XXXXX 00, 0000; accepted XXXXX 00, 0000.\\
This work was partially supported by NOVA LINCS (UIDB/04516/2020)
with the financial support of FCT I.P. and Project “Artificial Intelligence Fights Space Debris” No C626449889-0046305 co-funded by Recovery and Resilience Plan and NextGeneration EU Funds, www.recuperarportugal.gov.pt.
\includegraphics[scale=0.1]{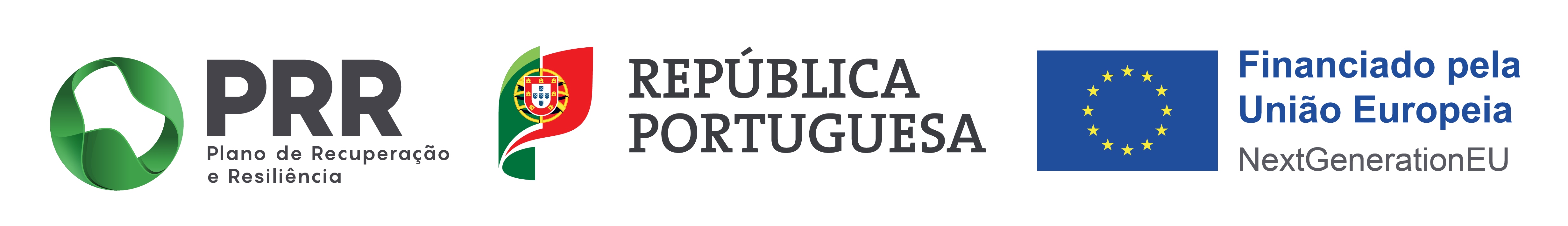}\\ }

\corresp{{\itshape (Corresponding author: Ricardo Ferreira)}.}

\authoraddress{Ricardo Ferreira, Filipa Valdeira and Cl\'{a}udia Soares are with NOVA School of Science and Technology
(e-mail: \href{rjn.ferreira@campus.fct.unl.pt}{rjn.ferreira@campus.fct.unl.pt},~\href{f.valdeira@fct.unl.pt}{f.valdeira@fct.unl.pt},~\href{claudia.soares@fct.unl.pt}{claudia.soares@fct.unl.pt}). Marta Guimar\~{a}es is with Neuraspace (e-mail: \href{marta.guimaraes@neuraspace.com}{marta.guimaraes@neuraspace.com}).}


\markboth{FERREIRA ET AL.}{MAXIMUM LIKELIHOOD ESTIMATOR FOR INITIAL ORBIT DETERMINATION IN LOW-EARTH ORBIT}
\maketitle

\begin{abstract}
With the increase in the number of active satellites and space debris in orbit, the problem of initial orbit determination (IOD) becomes increasingly important, demanding a high accuracy. Over the years, different approaches have been presented such as filtering methods (for example, Extended Kalman Filter), differential algebra or solving Lambert's problem. In this work, we consider a setting of three monostatic radars, where all available measurements are taken approximately at the same instant. This follows a similar setting as trilateration, a state-of-the-art approach, where each radar is able to obtain a single measurement of range and range-rate. Differently, and due to advances in Multiple-Input Multiple-Output (MIMO) radars, we assume that each location is able to obtain a larger set of range, angle and Doppler shift measurements. Thus, our method can be understood as an extension of trilateration leveraging more recent technology and incorporating additional data. We formulate the problem as a Maximum Likelihood Estimator (MLE), which for some number of observations is asymptotically unbiased and asymptotically efficient. Through numerical experiments, we demonstrate that our method attains the same accuracy as the trilateration method for the same number of measurements and offers an alternative and generalization, returning a more accurate estimation of the satellite's state vector, as the number of available measurements increases.
\end{abstract}



\section{Introduction}

Initial orbit determination (IOD) is a crucial step in the Space Situational Awareness (SSA) pipeline. The goal of IOD is to obtain a six-parameter state vector, at a particular instant, that is able to describe an orbit (e.g., position and velocity in the Cartesian reference frame or the Keplerian elements)~\cite{schutz2004statistical}. With the exponential increase in the number of Resident Space Objects (RSOs)~\cite{esa-numbers} and the crucial role that active satellites play in contemporary society, it is important to obtain the object's state vector in an accurate and timely manner.

When an orbital state is available (usually for cataloged objects that are constantly monitored), new observations can be used in a filtering approach, such as the Extended Kalman Filter (EKF)~\cite{schutz2004statistical,vallado2001fundamentals,lam2010analysis,pardal2011robustness}. With an initial orbital state, we can propagate a reference trajectory and through a first-order Taylor Expansion at each observable instant, we can formulate the problem as a linear system and apply the standard Kalman filter, which iteratively converges to the true state vector, given an accurate initialization~\cite{krener2003convergence}. Furthermore, the propagation of the reference orbit relies on a simplified model that does not take into account relevant factors, such as the behavior of perturbative forces (for example, atmospheric drag and solar radiation pressure). This introduces a structural uncertainty that is not taken into account~\cite{poore2016covariance}, as well as the integration of Ordinary Differential Equations (ODEs), which can be time-consuming, depending on the desired accuracy.

Nevertheless, the problem of initial orbit determination becomes especially relevant in cases where we observe an object for the first time, with no prior information about the object's orbital state.

Over the years, different approaches have been presented to the problem of initial orbit determination of near-Earth orbiting satellites, when there is no \textit{a priori} knowledge about the orbit of the object. Some of these approaches only consider angle observations to obtain three position vectors at different instants, such as Gauss's method and Double-r iteration~\cite{vallado2001fundamentals}. The accuracy of Gauss's method is very sensitive to the separation between observations, preferably less than 10 degrees. Double-r iteration, incorporated in a method proposed by Escobal~\cite{escobal1970methods} is able to handle observations that are days apart. To obtain the full-state vector, one can resort to Gibbs method~\cite{gibbs1889determination}, which performs best for larger time lengths between position vectors, or the Herrick-Gibbs method~\cite{herrick1971astrodynamics}, which was developed for smaller time lengths when the vectors are almost parallel~\cite{kaushik2016statistical}. Given three position vectors, these methods are able to obtain the velocity vector at the middle point. Another solution to obtain the velocity vector is to solve Lambert's problem~\cite{izzo2015revisiting,gooding1988solution,gooding1990procedure,lancaster1969unified,gooding1996new}, which determines the Keplerian elements given two position vector and the time period between the two positions, often called \textit{time of flight}. More recently, Qu et al. showed that by combining Doppler measurements with Lambert's problem, the authors provide an analytic solution, substituting the standard iterative algorithms to solve Lambert's problem. Another approach in the same vein of EKF is a Least Squares solution~\cite{schutz2004statistical}. However, instead of using each new observation to update the orbit, all the observations are jointly taken to improve the initial orbital state. While this approach is commonly used to determine the orbits of non-cooperative LEO satellites, the drawbacks of the EKF extend to the Least Squares solution~\cite{deng2022non,deng2023non}.

Lately, differential algebra (DA) has been exploited to improve initial orbit determination. With differential algebra, we can relate the uncertainty in measurements and the orbital state through a truncated power series (TPS) from Taylor polynomials of some order $k$. From this, we can resort to automatic domain splitting~\cite{wittig2015propagation}, which restricts the domain along each TPS direction when the truncation error surpasses some tolerance $\varepsilon$. So, we are able to obtain a compact region of possible solutions, in which the true orbit lies. In the literature, we find the application of these DA-based approaches for initial orbit determination under different settings, in particular, considering only angle measurements~\cite{pirovano2020probabilistic}, range and angle measurements~\cite{armellin2018probabilistic} and Doppler-only radars~\cite{losacco2023robust}. While~\cite{pirovano2020probabilistic,armellin2018probabilistic,losacco2023robust} resort to non-perturbed propagation of the orbital state to obtain a solution, Foss{\`a} et al. recently proposed the corresponding extensions for perturbed dynamical systems~\cite{fossa2024perturbed}.

Also in recent years, the use of Machine Learning (ML) has been explored for the problem of initial orbit determination~\cite{caldas2024machine}. {\"O}reng{\"u}l et al. study the application of Artificial Neural Networks to retrieve the Keplerian elements from angular observations~\cite{orengul2023artificial}. Schwab et al. examine the performance of multivariate Gaussian Process regression, for angles-only observations, to accurately determine the orbital state and quantify the uncertainty associated with the estimation~\cite{schwab2022angles}. Lee et al. use an ensemble of different ML models for range-rate and slant-range data generation~\cite{lee2018machine}. Relying on tracking data from a real station to train the models, the authors are able to estimate measurements from two virtual ground stations, creating additional data to be used in classic IOD methods.

Other approaches consider the setting where a small number of observations are available (\textit{too-short-arcs})~\cite{shang2019VSA,qu2022VSA,zhang2019initial}. State-of-the-art approaches gather new information from high-order kinematic parameters which can be obtained from time derivatives of radar's echo phase, however, this is only possible for Low-Earth objects with stable attitude.

In this work, we consider the scenario where all available measurements are taken approximately at the same instant (one-shot). This scenario is particularly important for the initial orbit determination (IOD) of small debris. Other approaches consider the same setting, such as the \emph{trilateration} method~\cite{escobal1970methods,vallado2001fundamentals,hough2012precise}, which bypass the need for assumptions about the dynamical systems, while accurately determining the state vector and associated uncertainty. The trilateration approach requires, simultaneously, three range measurements and three range-rate measurements~\cite{hough2012precise}. Recently, considering simultaneous time delay and Doppler shift measurements,~Ferreira et al. formulated the problem of initial orbit determination as a Weighted Least Squares problem~\cite{ferreira2023one}, directly obtaining the state of the object (position and velocity)
and the associated covariance matrix from Fisher’s Information Matrix. Different from the trilateration that considers a set of three monostatic radars, this method considers a multistatic radar system with $N$ transmitters and $M$ receivers.  

For our approach, we consider the same radar setting as the trilateration, i.e., a set of three monostatic radars, and that all the measurements are taken approximately at the same instant. Instead of considering one measurement of range and one of range-rate per radar, due to advances in Multiple-Input Multiple-Output (MIMO) radars~\cite{martinez2017mimo,misiurewicz2019mimo,forsythe2004multiple}, we assume that each location is able to obtain a larger set of range, angle and Doppler shift measurements. With a larger set of data, assuming that measurement noise follows a known distribution, we can formulate the problem of initial orbit determination as a Maximum Likelihood Estimator (MLE), which for large sets of data and high signal-to-noise ratio, is asymptotically unbiased and efficient~\cite{kay1993fundamentals}. By relaxing the cost function of our Maximum Likelihood Estimator, we can apply a block coordinate gradient descent algorithm to iteratively solve two convex problems, provably improving on the original non-convex problem.

\subsection{Contributions}

In comparison with the state-of-the-art approaches for the problem of initial orbit determination, our approach retains the following advantages and differences:

\begin{itemize}
    \item We formulate the problem as a Maximum Likelihood Estimator (MLE), which for some number of observations is asymptotically unbiased and asymptotically efficient~\cite{kay1993fundamentals};
    \item We do not rely on the propagation of a reference trajectory, thus avoiding linear approximations and simplifications of the physical systems;
    \item We show that our approach attains the same accuracy as the trilateration for the same number of measurements and the error decreases as the number of measurements increases. Therefore, our approach presents as a generalization of the trilateration method for an arbitrary number of measurements.
\end{itemize}

\section{Problem Formulation}

In this section, we describe the problem we want to solve. As previously mentioned, the goal is to determine a six-parameter state vector that describes the motion of a satellite, specifically the position, $x \in {\mathbb{R}}^{3}$, and velocity, $v \in {\mathbb{R}}^{3}$, in the Cartesian reference frame.

We consider a set of $N$ monostatic radars (capable of transmitting and receiving its own signal), located at $t_{i} \in {\mathbb{R}}^{3}$, for $i = 1, \ldots, N$, all able to observe the satellite at the same instant. We assume that each radar $i$ is able to collect a set of different measurements:

\begin{itemize}
    \item Range measurements, $d_i = \|x - t_i\| + \epsilon_{d_i}$;
    \item Angle measurements, $ u_i = \frac{x - t_i}{\|x - t_i\|} + \epsilon_{u_i}$;
    \item Doppler shift measurements
    \begin{align*}
        f_i = \frac{2 f_{c,i}}{c} \left(\frac{x - t_i}{\|x - t_i\|}\right)^T v + \epsilon_{f_i}.
    \end{align*}
\end{itemize}
where $c$ denotes the speed of light and $f_{c,i}$ denotes the carrier frequency of the signal from radar $i$. Range and Doppler shift noise is modeled as a Gaussian distribution with zero mean and variance $\sigma^2_{d_i}$ and $\sigma^2_{f_i}$, respectively, so $\epsilon_{d_i} \sim \mathcal{N}\left(0, \sigma^2_{d_i}\right)$ and $\epsilon_{f_i} \sim \mathcal{N}\left(0, \sigma^2_{f_i}\right)$.

Angle measurements are represented as unit-norm vectors and can be retrieved from elevation and azimuth measurements. Some methods model elevation and azimuth noise following a Gaussian distribution~\cite{yanez2017novel,fossa2024perturbed}, however, we assume that angle measurements are obtained as unit-norm directional vectors, therefore we model angle noise as a von Mises-Fisher distribution with mean direction zero and concentration parameter $\kappa_i$, as in~\cite{valdeira2024maximum}, which can better represent the distribution of directional data, so $\epsilon_{u_i} \sim VMF\left( 0, \kappa_i \right)$.

We can still relate the concentration parameter $\kappa_i$ with an equivalent standard deviation $\sigma_{u_i}$ as~\cite{mardia2000directional}

\begin{equation}
    \label{eq:angle-std}
    \centering
    \begin{aligned}
        \sigma_{u_i} = \sqrt{-2 \ln\left( 1 - \frac{1}{2 \kappa_i} - \frac{1}{8 \kappa_i^2} - \frac{1}{8 \kappa_i^3} \right)}.
    \end{aligned}
\end{equation}

Thus, the problem we intend to solve is to estimate the position, $x$, and velocity, $v$, of a satellite given the dataset of measurements

\begin{equation}
    \label{eq:measurement-dataset}
    \centering
    \begin{aligned}
        {\mathcal{D}} =& \underbrace{\left\{ d_i \in {\mathbb{R}} : i = 1, \dots, N \right\}}_{\text{range measurements}} \ \bigcup \\
        & \underbrace{\left\{ u_i \in {\mathbb{R}}^{3} : i = 1, \dots, N \right\}}_{\text{angle measurements}} \ \bigcup \\
        & \underbrace{\left\{ f_i \in {\mathbb{R}} : i = 1, \dots, N \right\}}_{\text{Doppler shift measurements}}.
    \end{aligned}
\end{equation}
By assuming that the noise is independent and identically distributed, we can obtain the maximum likelihood estimator by solving the optimization problem

\begin{equation}
    \label{eq:mle-estimator-problem}
    \centering
    \begin{aligned}
    \underset{x, v}{\mbox{minimize}} \quad & f_{\text{range}}(x) + f_{\text{angle}}(x) + f_{\text{doppler}}(x,v)
    \end{aligned}
\end{equation}
such that the functions are defined as

\begin{equation}
    \label{eq:f-range-function}
    \centering
    \begin{aligned}
    f_{\text{range}}(x) = \sum_{i=1}^{N} \frac{1}{2 \sigma_{d_i}^2} \left( \|x - t_i\| - d_i \right)^2
    \end{aligned}
\end{equation}
\begin{equation}
    \label{eq:f-angle-function}
    \centering
    \begin{aligned}
    f_{\text{angle}}(x) = -\sum_{i=1}^{N} \kappa_i u_i^T \frac{x - t_i}{\|x - t_i\|}
    \end{aligned}
\end{equation}
\begin{equation}
    \label{eq:f-doppler-function}
    \centering
    \begin{aligned}
    f_{\text{doppler}}(x,v) = \sum_{i=1}^{N} \frac{1}{2 \sigma_{f_i}^2} \left( \frac{2 f_{c,i}}{c} \left(\frac{x - t_i}{\|x - t_i\|}\right)^T v - f_i \right)^2,
    \end{aligned}
\end{equation}
for the range, angle and Doppler shift measurements, respectively. Problem~\eqref{eq:mle-estimator-problem} is non-convex, due to the range, angle and Doppler terms. For the range term, nonconvexity emerges from the square of distances smaller than $d_i$, i.e., when $\|x - t_i\| < d_i$. For the angle term, nonconvexity arises from the nonlinear denominator $\|x - t_i\|$, which is similar to the Doppler term that also involves a bilinear term with the two optimization variables $x$ and $v$.

We follow the relaxations proposed in~\cite{soares2015simple,soares2020range}, by approximating the functions $f_{\text{range}}$ and $f_{\text{angle}}$ with convex functions.

First, we introduce a set of new variables $y_i \in {\mathbb{R}}^{3}$ for $i = 1, \ldots, N$ and then rewrite the terms of the function $f_{\text{range}}$ to an equivalent formulation

\begin{equation}
    \label{eq:relaxed-f-range-function}
    \centering
    \begin{aligned}
    \left( \|x - t_i\| - d_i \right)^2 = \inf_{\|y_i\| = d_i} \left\| x - t_i - y_i \right\|^2.
    \end{aligned}
\end{equation}
To obtain a convex problem, we relax the constraint as $\|y_i\| \leq d_i$, as in~\cite{soares2015simple}. We can use the new set of variable $y_i$ and the range measurements to relax the angular terms, as in~\cite{soares2020range}. In this way, we can replace the terms $\left(x - t_i\right)$ by the variable $y_i$ and the denominator $\left\|x - t_i\right\|$, which represents the distance between the variable $x$ and the radar located at $t_i$, by the corresponding range measurement, $d_i$. Thus, the new formulation is given by

\begin{equation}
    \label{eq:relaxed-f-angle-function}
    \centering
    \begin{aligned}
    \hat{f}_{\text{angle}}(y) = -\sum_{i=1}^{N} \frac{\kappa_i}{d_i} u_i^T y_i.
    \end{aligned}
\end{equation}
where $y$ can be seen as the concatenation of all the new variables $y_i$, i.e., $y \coloneqq (y_1, \ldots, y_N)$. The same idea can be applied to the term $f_{\text{doppler}}$ leading to

\begin{equation}
    \label{eq:relaxed-f-doppler-function}
    \centering
    \begin{aligned}
    \hat{f}_{\text{doppler}}(y,v) = \sum_{i=1}^{N} \frac{\beta_i^2}{2} \left( \omega_i y_i^T v - f_i \right)^2,
    \end{aligned}
\end{equation}
such that $\beta_i = \frac{1}{\sigma_{f_i}}$ and $\omega_i = \frac{2 f_{c,i}}{c d_i}$. So, the new problem is written as

\begin{equation}
    \label{eq:relaxed-mle-estimator-problem}
    \centering
    \begin{aligned}
    &\underset{x, y, v}{\mbox{minimize}} \quad \hat{f}_{\text{range}}(x,y) + \hat{f}_{\text{angle}}(y) + \hat{f}_{\text{doppler}}(y,v) \\
    &\mbox{subject to } \quad \|y_i\| \leq d_i, \quad \mbox{for } i = 1, \ldots, N,
    \end{aligned}
\end{equation}
where 

\begin{equation}
    \label{eq:relaxed-f-range-function-2}
    \centering
    \begin{aligned}
    \hat{f}_{\text{range}}(x,y) = \sum_{i=1}^{N} \frac{\alpha_i^2}{2} \left\|x - t_i - y_i \right\|^2,
    \end{aligned}
\end{equation}
such that $\alpha = \frac{1}{\sigma_{d_i}}$. Problem~\eqref{eq:relaxed-mle-estimator-problem} is still non-convex due to the term $\hat{f}_{\text{doppler}}$, because of the bilinear terms $y_i^T v$. However, we can sidestep this problem by solving it through a block coordinate descent approach~\cite{tseng2001convergence,luo1993error,beck2013convergence}. Iteratively, we minimize our problem for the variables $y_i$ for $i = 1, \ldots, N$, and then we apply the same process but for the variables $x$ and $v$.

\section{Block Coordinate Descent}

To solve our problem through a block coordinate descent approach~\cite{peng2016coordinate,wright2015coordinate,shi2016primer}, we consider two major blocks: we simultaneously minimize for $x$ and $v$, and posteriorly minimize for the variables $y_i$, for $i = 1, \ldots, N$.

\subsection{Minimization over $(x, v)$}

In order to minimize the variables $x$ and $v$, we analyze the unconstrained convex problem

\begin{equation}
    \label{eq:x,v-general-problem-formulation}
    \begin{aligned}
        &\underset{x, v}{\mbox{minimize}} \quad F_y(x, v) = \hat{f}_{\text{range}}(x,y) + \hat{f}_{\text{doppler}}(y,v) \\
    \end{aligned}
\end{equation}
where we fix the variables $y_i$, for $i = 1, \ldots, N$, to their current value at the corresponding iteration. This problem is a Quadratic Program (QP)~\cite{boyd2004convex} and has a closed-form solution, which can be obtained by equating the gradient of the cost function $F_y(x, v)$ to zero. 

\begin{restatable}{prop}{propone}
\label{prop:first-proposition}
    The solution to problem~\eqref{eq:x,v-general-problem-formulation} is given by $x^{*} = W^{-1} q$ and $v^{*} = \left( Y^T B Y \right)^{-1} Y^T B f$, such that

    \begin{align*}
        W = \sum_{i=1}^{N} \alpha_{i}^2 I, \quad q = \sum_{i=1}^{N} \alpha_{i}^2 (t_i + y_i), \\
    \end{align*}
    \begin{align*}
        Y = \begin{bmatrix}
            \omega_1 y_1^T \\
            \vdots \\
            \omega_N y_N^T
        \end{bmatrix}, \quad
        B = \begin{bmatrix}
            \beta_1^2 & & \\
            & \ddots & \\
            & & \beta_N^2
        \end{bmatrix}, \quad
        f = \begin{bmatrix}
            f_1 \\
            \vdots \\
            f_N
        \end{bmatrix}.
    \end{align*}
\end{restatable}

\paragraph*{Proof.} The terms $\hat{f}_{\text{range}}$ and $\hat{f}_{\text{doppler}}$ can be written as

\begin{equation}
    \label{eq:reformulated-f-range-and-doppler-function}
    \centering
    \begin{aligned}
    \hat{f}_{\text{range}}(x,y) &= \frac{1}{2} \sum_{i=1}^{N} \left\|\alpha_i x - \alpha_i (t_i + y_i)\right\|^2 \\
    \hat{f}_{\text{doppler}}(y,v) &= \frac{1}{2} \left(Y v - f\right)^T B \left(Y v - f\right)
    \end{aligned}
\end{equation}
We then obtain the partial derivatives of $F_y$ as

\begin{equation}
    \label{eq:x,v-gradient}
    \begin{aligned}
        \nabla_x F_y(x, v) &= \sum_{i=1}^{N} \alpha_i I \left( \alpha_i x - \alpha_i (t_i + y_i) \right) \\
        &= W x - q \\
        \nabla_v F_y(x, v) &= Y^T B \left( Y v - f \right)
    \end{aligned}
\end{equation}
By equating both gradients to zero, we obtain the optimal solution as

\begin{align*}
    x^{*} = W^{-1} q, \quad
    v^{*} = \left( Y^T B Y \right)^{-1} Y^T B f. \qed
\end{align*}

\subsection{Minimization over $y$}

To minimize for the variables $y_i$, for $i = 1, \ldots, N$, our problem becomes

\begin{equation}
    \label{eq:y_i-problem-formulation}
    \begin{aligned}
        &\underset{y}{\mbox{minimize}} \quad F_{x,v}(y) = \hat{f}_{\text{range}}(x,y) + \hat{f}_{\text{angle}}(y) + \hat{f}_{\text{doppler}}(y,v)  \\
        &\mbox{subject to} \quad \|y_i\|^2 \leq d_i^2 \quad \mbox{for } i = 1, \ldots, N  \\
    \end{aligned}
\end{equation}
where we fix the values for $x$ and $v$ and solve for the variables $y_i$, for $i = 1, \ldots, N$. However, we see that the problem is separable for each variable $y_i$

\begin{equation}
    \label{eq:y_i-problem-formulation-separable}
    \begin{aligned}
        F_{x,v}(y) &= \hat{f}_{\text{range}}(x,y) + \hat{f}_{\text{angle}}(y) + \hat{f}_{\text{doppler}}(y,v)  \\
        &= \sum_{i=1}^{N} \underbrace{ \frac{\alpha_i^2}{2} \left\|a_i - y_i \right\|^2 - b_i^T y_i + \frac{\beta_i^2}{2} \left( \omega_i y_i^T v - f_i \right)^2 }_{g_i\left(y_i\right)},
    \end{aligned}
\end{equation}
such that $a_i = x - t_i$ and $b_i = \frac{\kappa_i}{d_i} u_i$, for $i = 1, \ldots, N$. So, we can independently solve for each variable $y_i$, by solving a Quadratically Constrained Quadratic Program (QCQP)~\cite{boyd2004convex}, with $g_i(y_i)$ as the objective function. We can reformulate our objective function as a quadratic function in canonical form as

\begin{equation}
    \label{eq:y_i-quadratic-objective}
    \begin{aligned}
        g_i\left( y_i \right) = \frac{1}{2} y_i^T A_i y_i + p_i^T y_i + \zeta_i
    \end{aligned}
\end{equation}
where

\begin{equation}
    \label{eq:y_i-quadratic-objective-aux-variables}
    \begin{aligned}
        A_i &= \alpha_i^2 I + \omega_i^2 \beta_i^2 v v^T \\
        p_i &= b_i - \alpha_i^2 a_i - \beta_i^2 \omega_i f_i v \\
        \zeta_i &= \frac{1}{2} \left(\alpha_i^2 a_i^T a_i + \beta_i^2 f_i^2\right) \\
    \end{aligned}
\end{equation}
So, our convex optimization problem becomes

\begin{equation}
    \label{eq:y_i-problem-formulation-trs}
    \begin{aligned}
        &\underset{y_i}{\mbox{minimize}} \quad \frac{1}{2} y_i^T A_i y_i + p_i^T y_i  \\
        &\mbox{subject to} \quad \|y_i\| \leq d_i  \\
    \end{aligned}
\end{equation}
where we discard the term $\zeta_i$, since it is a constant. This problem is convex because $A_i$ is positive semidefinite. For all $x \in {\mathbb{R}}^{3}$, we have

\begin{align*}
    x^T A_i x = \alpha_{i}^2 \|x\|^2 + \omega_i^2 \beta_i^2 (x^T v)^2 \geq 0,
\end{align*}
since $\alpha_i, \omega_i, \beta_i$ are always positive.

More than a QCQP, Problem~\eqref{eq:y_i-problem-formulation} is part of a class of problems called Trust Region Subproblem (TRS)~\cite{yuan2000review,adachi2017solving,fortin2004trust}, whose algorithms (denominated \textit{trust region methods}) are known for their strong convergence properties.

For the Trust Region Subproblem (TRS), the primal and dual optimal points, $y_i^{*}$ and $\lambda_i^{*}$, respectively, satisfy the optimality conditions~\cite{fortin2004trust,sorensen1982newton}

\begin{subequations}
    \label{eq:y_i-problem-optimality-conditions}
    \begin{align}
        \label{eq:primal_feasibility-1}
        &\left\| y_i^{*} \right\| \leq d_i \\
        \label{eq:dual_feasibility-1}
        &\left(A_i + \lambda_i^{*} I\right) y_i^{*} = -p_i \\
        \label{eq:dual_feasibility-2}
        &A_i + \lambda_i^{*} I \succeq 0 \\
        \label{eq:dual_feasibility-3}
        &\lambda_i^{*} \geq 0 \\
        \label{eq:complementary-slackness-1}
        &\lambda_i^{*} \left( \left\| y_i^{*} \right\|^2 - d_i^2 \right) = 0
    \end{align}
\end{subequations}
From ~\eqref{eq:dual_feasibility-1}, we have that the primal solution is given by

\begin{equation}
    \label{eq:y_i-primal-solution}
    \begin{aligned}
        y_i^{*} = -\left(A_i + \lambda_i^{*} I\right)^{-1} p_i ,
    \end{aligned}
\end{equation}
so we must analyze in which cases $\left(A_i + \lambda_i^{*} I\right)$ is invertible. 

\begin{restatable}{prop}{proptwo}
\label{prop:second-proposition}
    The matrix $\left(A_i + \lambda_i^{*} I\right)$ is invertible for all $\lambda_i^{*} \geq 0$.
\end{restatable}
\noindent A proof of this proposition can be found in the appendix. So, now we have two cases to satisfy the first and last two optimal conditions (\eqref{eq:primal_feasibility-1},~\eqref{eq:dual_feasibility-3} and~\eqref{eq:complementary-slackness-1}):

\begin{itemize}
    \item[(1)] $\lambda_i^{*} = 0$ and $\left\| y_i^{*} \right\| = \left\|-A_i^{-1} p_i \right\| < d_i$;
    \item[(2)] $\lambda_i^{*} > 0$ and $\left\| y_i^{*} \right\| = \left\| -\left( A_i + \lambda_i^{*} I \right)^{-1} p_i \right\| = d_i$.
\end{itemize}
The first case is easy to confirm. If we compute $\left\|-A_i^{-1} p_i \right\|$ and it is less or equal to $d_i$, then we have the primal and dual optimal solutions $\left(y_i^{*}, \lambda_i^{*}\right) = \left(-A_i^{-1} p_i, 0\right)$. Otherwise, we must find the positive $\lambda_i^{*}$ that satisfies the condition $\left\|-\left( A_i + \lambda_i^{*} I \right)^{-1} p_i \right\| = d_i$. 

In~\cite{adachi2017solving}, Adachi et al. jointly with the work of Gander et al.~\cite{gander1989constrained} proved that a solution $(y_{i}^{*}, \lambda_{i}^{*})$ of the TRS (Problem~\eqref{eq:y_i-problem-formulation-trs}) satisfying the optimality conditions~\eqref{eq:primal_feasibility-1}-~\eqref{eq:complementary-slackness-1}, for a primal solution on the boundary (i.e., $\|y_{i}^{*}\| = d_{i}$), the optimal dual solution $\lambda_{i}^{*}$ is equal to the largest real eigenvalue of the nonsymmetric eigenvalue problem

\begin{equation}
    \label{eq:optimal-dual-solution}
    \begin{aligned}
        \begin{bmatrix}
            -A_i & I \\
            \frac{1}{d_{i}^2} p p^T & -A_i
        \end{bmatrix} z = \lambda z
    \end{aligned}
\end{equation}
To compute the primal solution, $y_{i}^{*}$, we start by assessing if the largest real eigenvalue of Problem~\eqref{eq:optimal-dual-solution} is positive. If true, then the primal solution is $y_i^{*} = \left( A_i + \lambda_i^{*} I \right)^{-1} p_i$. Otherwise, the optimal value for $\lambda_i^{*}$ is zero and the primal solution is given by $y_i^{*} = A_i^{-1} p_i$. The entire algorithm is made explicit in Algorithm~\ref{alg:block-coordinate-descent}.

\begin{algorithm}
\caption{Block Coordinate Descent}\label{alg:block-coordinate-descent}
\textbf{Input:} Dataset $\mathcal{D} = \left(d_1, \ldots, d_N, u_1, \ldots, u_N, f_1, \dots, f_N\right)$ 
\begin{algorithmic}
\State Initialize $x, v, y_1, \ldots, y_N$
\While{stopping criterion is not met}
    \For{$i = 1, \ldots, N$}
        \State Find largest real eigenvalue $\lambda_i^{*}$ of the eigenvalue problem $\begin{bmatrix}
            -A_i & I \\
            \frac{1}{d_{i}^2} p p^T & -A_i
        \end{bmatrix} z = \lambda z$ \Comment{Problem.~\eqref{eq:optimal-dual-solution}}
        \If{$\lambda_i > 0$}
            \State $y_i = \left( A_i + \lambda_i^{*} I \right)^{-1} p_i$
        \Else
            \State $y_i = A_i^{-1} p_i$
        \EndIf
    \EndFor
    
    \State $x = W^{-1} q$ \Comment{Proposition~\ref{prop:first-proposition}}
    \State $v = \left( Y^T B Y \right)^{-1} Y^T B F$ \Comment{Proposition~\ref{prop:first-proposition}}
    
\EndWhile

\State\Return $\left\{x, v, y_1, \ldots, y_N\right\}$
\end{algorithmic}
\end{algorithm}

\section{Numerical Results}

For our numerical experiments, we consider a simulated scenario with five different objects. In Table~\ref{tab:orbital-information}, we can find the orbital information of the space objects used in our numerical experiments. In Table~\ref{tab:orbital-information}, $a$ denotes the semi-major axis in kilometers, $e$ represents the eccentricity of the orbit, and $i$, $\Omega$ and $\omega$ denote the inclination, the right ascension of the ascending node and the argument of the periapsis (in degrees), respectively. The only element that is not described is the mean anomaly, which we consider to be zero for all objects.

\begin{table}[ht!]
    \centering
    \caption{Keplerian elements of the selected objects for numerical experiments.}
    \label{tab:orbital-information}
    \begin{tabular}{|c|c|c|c|c|c|c|}
    \hline
    Object & $a$ & $e$ & $i$ & $\Omega$ & $\omega$ \\
    \hline
    Object \#1 & 6913.9278 & 0.0106 & 97.1377 & 66.7240 & 79.0900 \\
    Object \#2 & 6886.5427 & 0.0003 & 97.4457 & 68.2327 & 72.8300 \\
    Object \#3 & 6886.5577 & 0.0002 & 97.4460 & 67.7949 & 74.2700 \\
    Object \#4 & 7151.1996 & 0.0020 & 95.9746 & 68.1057 & 77.8000 \\
    Object \#5 & 6860.4158 & 0.0076 & 93.9043 & 64.4680 & 75.0700 \\
    \hline
    \end{tabular}
\end{table}
To compare our approach with the trilateration method, we consider three distinct positions where the radars are placed:

\renewcommand{\arraystretch}{1.5}
\begin{table}[h!]
    \centering
    \caption{Geographical coordinates of the radars of the considered multi-static radar system, as well as the carrier frequencies of the transmitter antennas.}
    \label{tab:lat-lon-table}
    \begin{tabular}{|c|c|c|c|}
    \hline
    & Latitude & Longitude & $f_c$ \\
    \hline
    $t_1$ & 72.986276º & 40.916634º & $1215$ MHz \\
    \hline
    $t_2$ & 74.986276º & 48.916634º & $1280$ MHz \\
    \hline
    $t_3$ & 75.986276º & 38.916634º & $1333$ MHz \\
    \hline
    \end{tabular}
\end{table}
These locations were chosen such that the radars have a clear vision of all different objects described in Table~\ref{tab:orbital-information}. From these locations, we assume that each location is able to obtain between one and five measurements of each type (range, angle, and Doppler shift), which can therefore be seen as a set between three and 15 radars, each one acquiring a measurement of each type. The first case would be comparable to the trilateration approach, where we are able to obtain one measurement of range and one of range-rate for each location (in addition, we also have a measurement of the angle). The latter would translate into obtaining a set of 45 measurements since each location can acquire 5 measurements of each type, which can be seen as a group of 15 radars.

So, for each set of radars (3, 6, 9, 12, and 15), we run $S = 100$ simulations and we compute the error for the position and velocity estimation as

\begin{equation}
    \centering
    \begin{aligned}
    \varepsilon_{n}^{i} = \norm{\Hat{x}_{n}^{i} - x_n^0}^2, \quad &\mbox{for } n = 1, \ldots, N, \\
    &\mbox{and } i = 1, \ldots, S,
    \end{aligned}
\end{equation}
where $N$ is the number of objects (which in our case $N = 5$), $\Hat{x}_{n}^{i}$ denotes the estimated position or velocity of the object $n$ for the $i$-th simulation, and $x_n^0$ denotes its true value. We simulate noisy range, angle and Doppler shift measurements. We consider the angle measurements' noise to follow a von Mises-Fisher distribution with mean direction zero and concentration parameter $\kappa_i = 10^{9}$ (equivalently, $\sigma_{u_i} = 0.0018^{o}$). We consider that the noise for range and Doppler shift measurements have mean zero and standard deviation $\sigma_{d_i} = 10^{-1} \, m$ and $\sigma_{f_i} = 10 \, Hz$, respectively. The noise for range measurements was chosen as a conservative noise of satellite laser ranging (SLR) observations~\cite{sosnica2015satellite}. The noise for Doppler shift measurements was based on the literature for orbit determination using Doppler shift measurements~\cite{deng2022non}.

To evaluate the robustness of our method, we consider three different scenarios where the noise for range and Doppler shift measurements follow: (1) a Gaussian distribution, (2) a Cauchy distribution and (3) a Laplace distribution.

\begin{figure}[htb]
    \centering 
\begin{minipage}[t]{.48\textwidth}
\begin{subfigure}{\textwidth}
  \includegraphics[width=\linewidth]{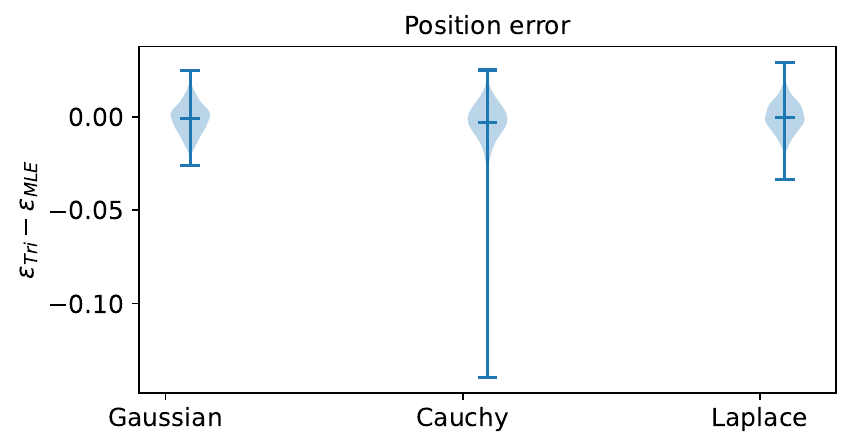}
  \label{fig:3.1}
\end{subfigure}\hfil 
\begin{subfigure}{\textwidth}
  \includegraphics[width=\linewidth]{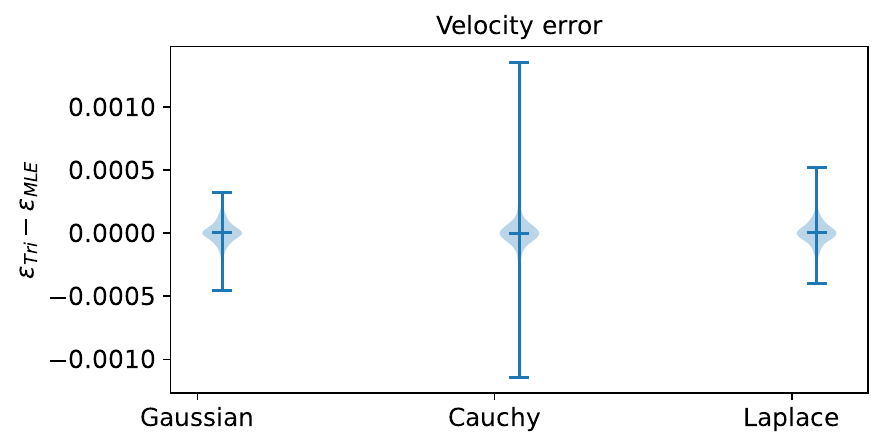}
  \label{fig:4.1}
\end{subfigure}\hfil 
\end{minipage}
\caption{Violin plot of the difference between the estimation errors of the trilateration method, $\varepsilon_{Tri}$, and the estimation errors of our proposed approach $\varepsilon_{MLE}$ considering three radars and measurements corrupted different noise distributions: Gaussian, Cauchy, and Laplace. For all types of noise distribution, we see that the estimation of both methods is similar as the difference between the errors is very close to zero, for both position and velocity.}
\label{fig:rmse-tri-vs-mle}
\end{figure}
In Fig.~\ref{fig:rmse-tri-vs-mle}, we consider the same setting for both trilateration and our approach, i.e., the same number of radars ($N=3$). For different types of noise (Gaussian, Cauchy and Laplace) we can see the violin plots of the difference between the estimation errors of the trilateration method, $\varepsilon_{Tri}$, and the estimation errors of our proposed approach, $\varepsilon_{MLE}$. These results show that our approach is equivalent to the trilateration method under these three types of noise since the difference between the errors is very close to zero.

\begin{figure*}
    \includegraphics[width=0.48\linewidth]{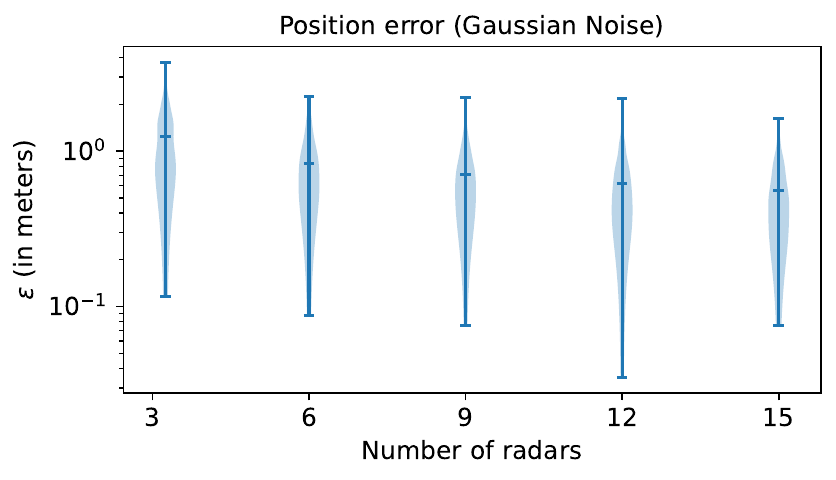} 
    \includegraphics[width=0.48\linewidth]{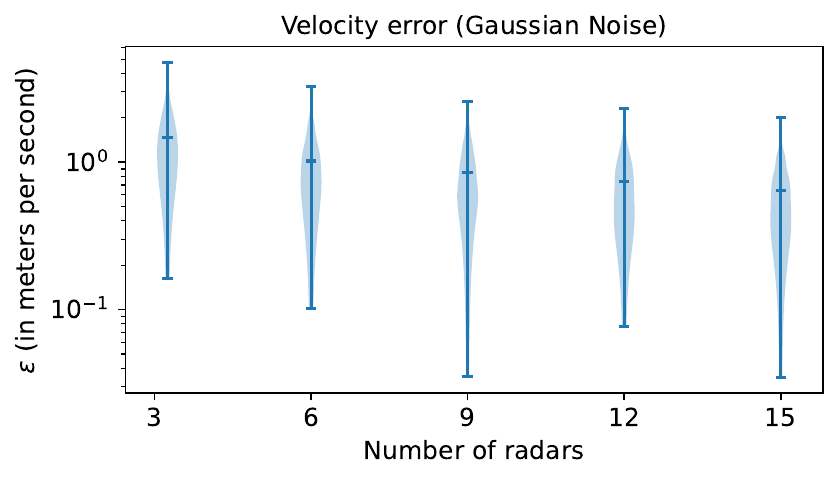}\par
    \includegraphics[width=0.48\linewidth]{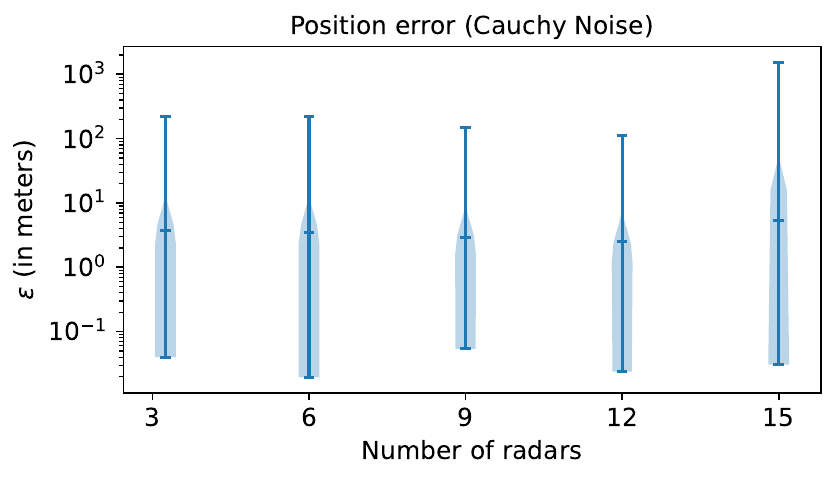} 
    \includegraphics[width=0.48\linewidth]{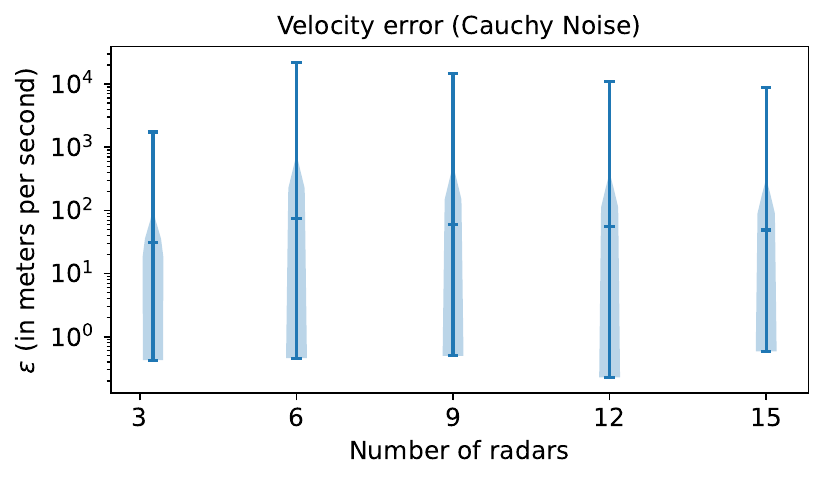}\par
    \includegraphics[width=0.48\linewidth]{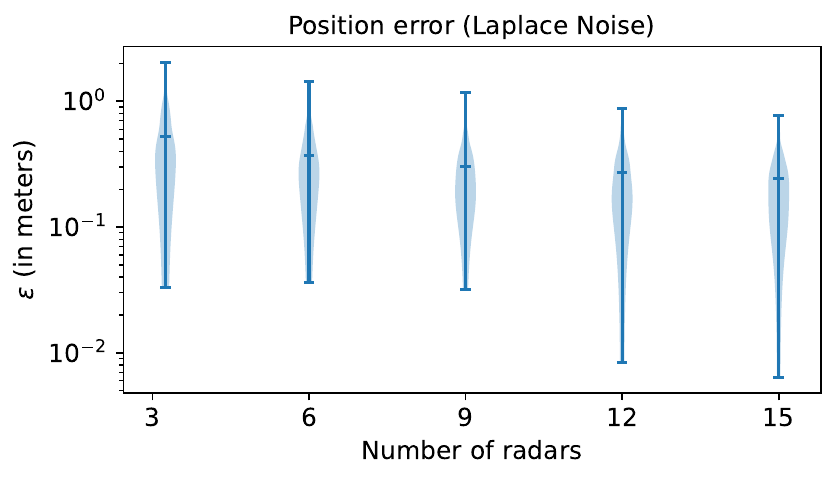} 
    \includegraphics[width=0.48\linewidth]{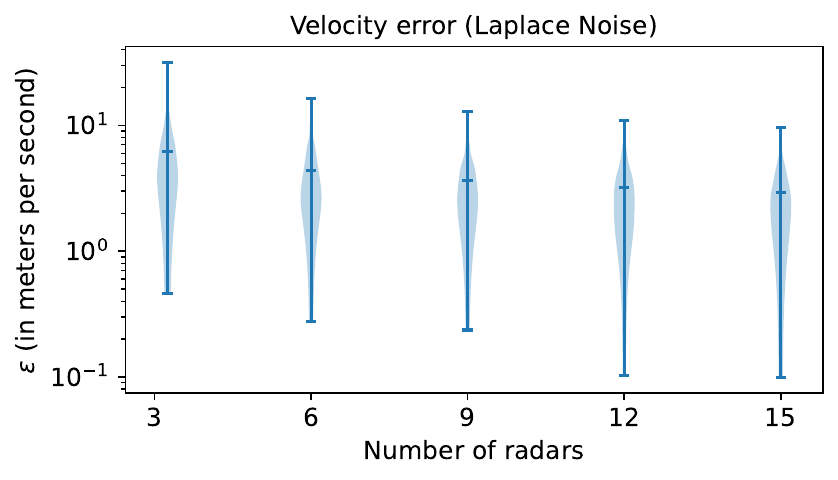}\par
\caption{Violin plot of the errors $\varepsilon$ of state estimate for Gaussian noise (in the first row), Cauchy Noise (in the second row), and Laplace Noise (in the third row). We see that for both Gaussian and Laplace noise, as the number of radars increases (therefore, measurements), the approximate Maximum Likelihood Estimator (MLE) becomes more accurate, as the mean of the errors $\varepsilon$ tends to decrease for both position and velocity estimation. Also, we can see that the uncertainty tends to decrease as the number of measurements increases. When measurements are corrupted by Cauchy noise, the error does not decrease and more outliers are generated.}
\label{fig:rmse-types-of-noise}
\end{figure*}
In Fig.~\ref{fig:rmse-types-of-noise}, we can see the evolution of the error as the number of radars (and consequently, measurements) increases, for different types of noise. As the number of radars increases, our estimator attains a lower error between our estimation and the true state (for both position and velocity), for Gaussian and Laplace noise. When measurements are corrupted by Laplace noise, the estimation error of the velocity seems to be higher than the estimation error of the position. This is different from the case of Gaussian noise, where the position and velocity estimations achieve the same level of accuracy.

For the case where the noise follows a Cauchy distribution, more outliers are generated and our method does not demonstrate robustness to deal with measurements corrupted by Cauchy noise. From these results, we conclude that our approach offers an alternative and generalization to the trilateration method, returning a more accurate estimation of the satellite's state vector, as the number of available measurements increases, demonstrating robustness for Gaussian and Laplace noise. Also, we can see that the uncertainty, for Gaussian and Laplace noise, tends to decrease as the number of measurements increases.

Our method attains the same accuracy for the same setting as the trilateration method, but in contrast is prepared to accept a larger number of measurements, demonstrating a smaller estimation error and uncertainty as the number of measurements increases.

\subsection{Robustness to different noise levels}

The second aspect we tested is the robustness of our approach to different noise levels. For all experiments, we consider the same setting for both trilateration and our approach, i.e., the same number of radars ($N=3$), and just alter the values of the noise applied to the measurements. As before, we run $S = 100$ simulations.

\begin{figure}[htb]
    \centering 
\begin{minipage}[t]{.48\textwidth}
\begin{subfigure}{\textwidth}
  \includegraphics[width=\linewidth]{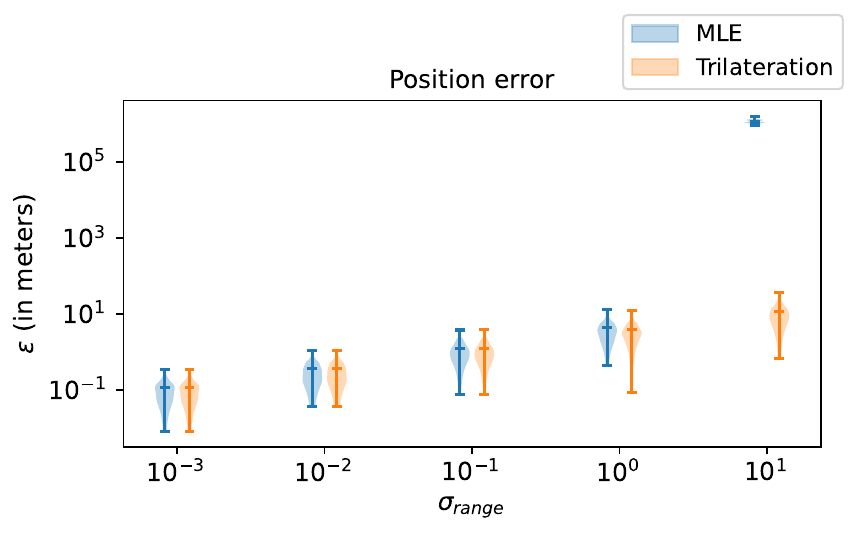}
  \label{fig:7}
\end{subfigure}\hfil 
\begin{subfigure}{\textwidth}
  \includegraphics[width=\linewidth]{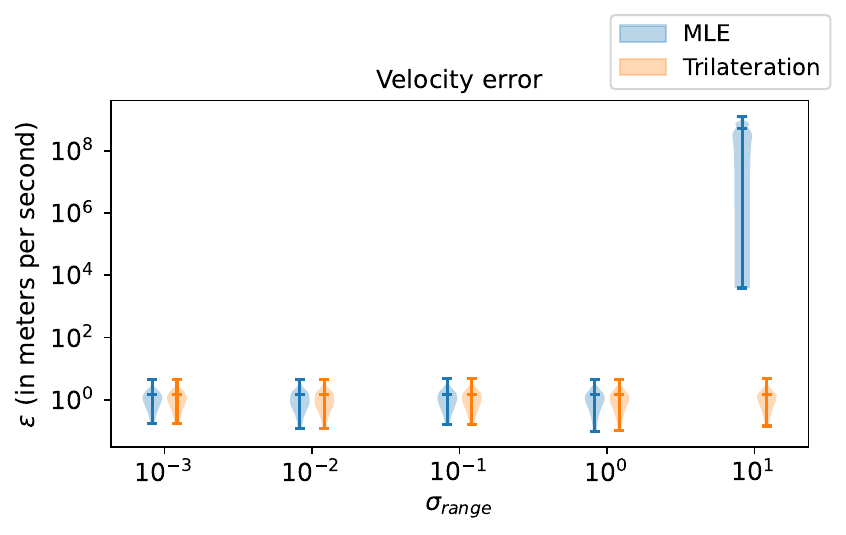}
  \label{fig:8}
\end{subfigure}\hfil 
\end{minipage}
\caption{Violin plot of the errors $\varepsilon$ of state estimate for different noise levels for the range measurements. As the range noise increases, we see that the position error increases sublinearly and the velocity error remains constant, with our MLE approach attaining the same accuracy as the trilateration method for noise up to $1 \,$ m.}
\label{fig:rmse-gaussian-range-noise-levels}
\end{figure}
\paragraph*{Range measurements} In the first experiment, we fixed $\kappa_i = 10^{9}$ and $\sigma_{f_i} = 10 \, Hz$, and computed the errors $\varepsilon$ for the position and velocity estimation of both our method and the trilateration approach, for the levels of noise $\sigma_{d_i} \in \{10^{-3}, 10^{-2}, 10^{-1}, 10^{0}, 10^{1}\} \, m$.

In Fig.~\ref{fig:rmse-gaussian-range-noise-levels}, we can observe that our estimation for both position and velocity attains the same accuracy as the trilateration approach for range noise up to $1 \, m$. The position error increases sublinearly, while the velocity error remains constant. This behaviour is expected since the range measurement directly impacts the position estimation, while the velocity estimation has a higher ``dependence'' on the Doppler shift measurements.

\begin{figure}[htb]
    \centering 
\begin{minipage}[t]{.48\textwidth}
\begin{subfigure}{\textwidth}
  \includegraphics[width=\linewidth]{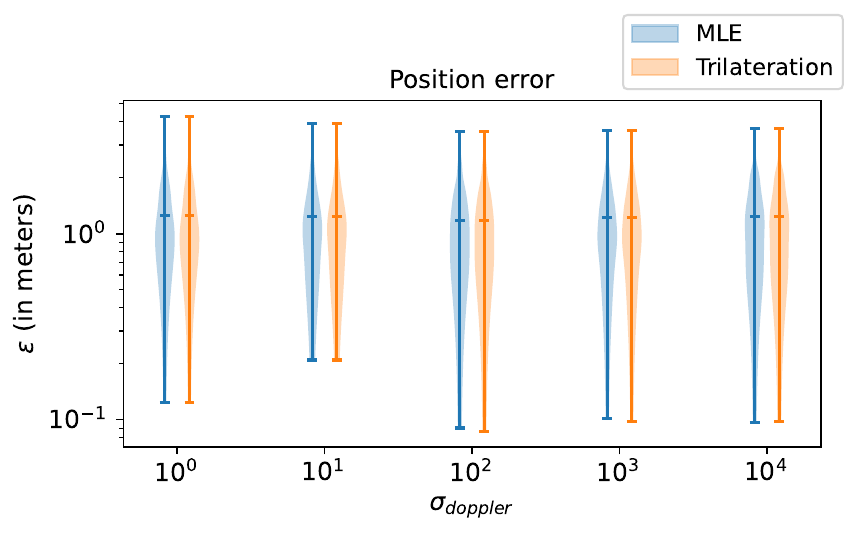}
  \label{fig:9}
\end{subfigure}\hfil 
\begin{subfigure}{\textwidth}
  \includegraphics[width=\linewidth]{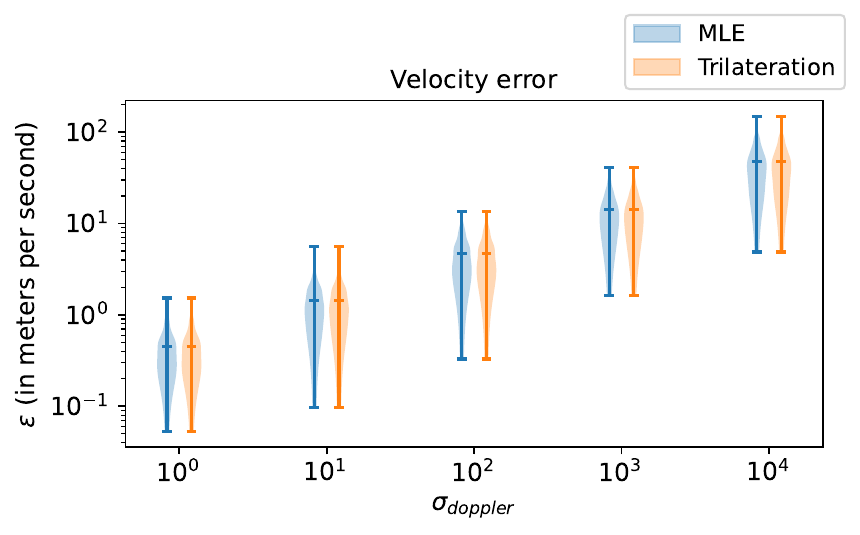}
  \label{fig:10}
\end{subfigure}\hfil 
\end{minipage}
\caption{Violin plot of the errors $\varepsilon$ of state estimate for different noise levels for the Doppler shift measurements. As the Doppler shift noise increases, we see that the position error remains constant and the velocity error increases sublinearly, with our MLE approach presenting a great level of robustness for the different noise levels, attaining the same accuracy as the trilateration method.}
\label{fig:rmse-gaussian-doppler-noise-levels}
\end{figure}
\paragraph*{Doppler shift measurements} For the second experiment, we fixed $\kappa_i = 10^{9}$ and $\sigma_{d_i} = 10^{-1} \, m$, and computed the errors $\varepsilon$ for the position and velocity estimation of both our method and the trilateration approach, for the levels of noise $\sigma_{f_i} \in \{10^{0}, 10^{1}, 10^{2}, 10^{3}, 10^{4}\} \, Hz$.

In Fig.~\ref{fig:rmse-gaussian-doppler-noise-levels}, we can observe the ``dependence'' of the velocity estimation on the Doppler shift measurements. While the position error remains constant for different levels of Doppler shift measurements noise, we see that the velocity error sublinearly increases as the noise increases. Also, as before, our approximate Maximum Likelihood Estimator achieves the same accuracy as the trilateration method for different levels of noise.

\begin{figure}[htb]
    \centering 
\begin{minipage}[t]{.48\textwidth}
\begin{subfigure}{\textwidth}
  \includegraphics[width=\linewidth]{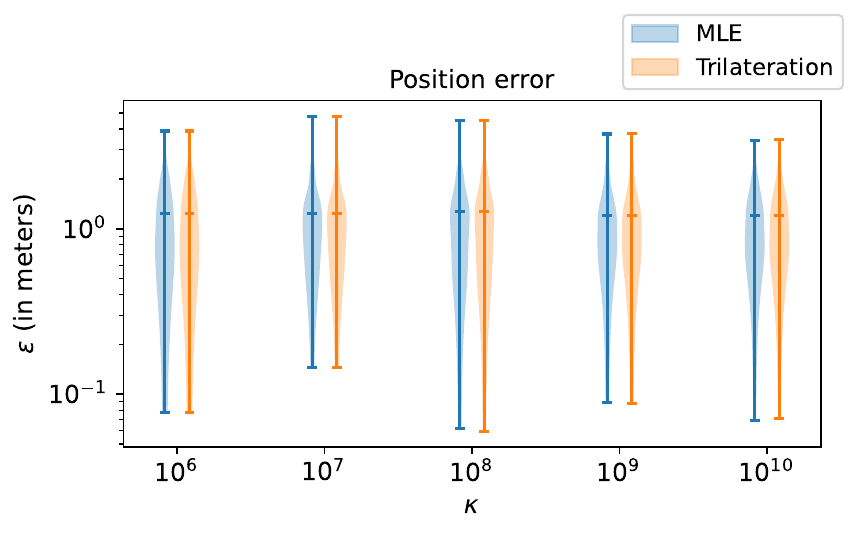}
  \label{fig:11}
\end{subfigure}\hfil 
\begin{subfigure}{\textwidth}
  \includegraphics[width=\linewidth]{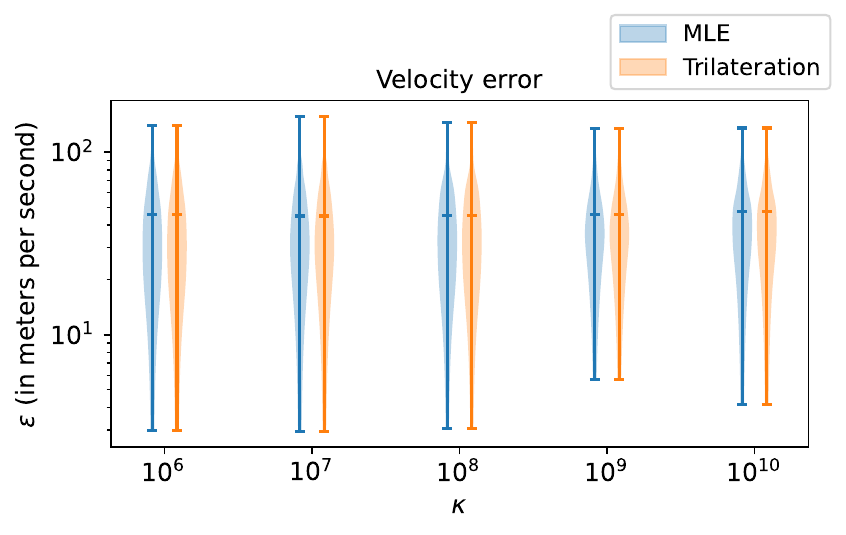}
  \label{fig:12}
\end{subfigure}\hfil 
\end{minipage}
\caption{Violin plot of the errors $\varepsilon$ of state estimate for different noise levels for the Doppler shift measurements. As $\kappa$ increases, we see that both position and velocity errors remain constant, with our MLE approach attaining the same accuracy as the trilateration method for the different noise levels.}
\label{fig:rmse-gaussian-bearing-noise-levels}
\end{figure}
\paragraph*{Angle measurements} Lastly, we fixed $\sigma_{d_i} = 10^{-1} \, m$ and $\sigma_{f_i} = 10 \, Hz$, and compute the errors $\varepsilon$ for the position and velocity estimation of both our method and the trilateration approach, for the levels of noise $\kappa \in \{10^{6}, 10^{7}, 10^{8}, 10^{9}, 10^{10}\}$.

Finally, in Fig.~\ref{fig:rmse-gaussian-bearing-noise-levels}, we can see that both position and velocity error remain constant for the different values for $\kappa$. This experiment shows that our algorithm is robust for a concentration parameter between $10^{6}$ and $10^{10}$, achieving the same accuracy as the trilateration method.

These three experiments show the equivalence between our approach and the trilateration method and the level of robustness for different levels of noise. For the same number of measurements, we show for different levels of noise in range, angle, and Doppler measurements, that our proposed approach can attain the same accuracy as the trilateration.


\section{Conclusions and Future Work}

With this work, we present an approximate Maximum Likelihood estimator for the problem of initial orbit determination for Low-Earth Orbit, which for a large set of data is asymptotically unbiased and asymptotically efficient.

Our estimator is able to provide the complete state vector of the object (position and velocity). Besides, our method does not need to propagate a reference trajectory, bypassing simplifications of the nonlinear dynamical system and avoiding time-consuming steps such as ODE integration.

We experimentally compare our approach with the trilateration method and with different numbers of measurements per radar. We show the equivalence between our proposed approach and the trilateration under three different types of noise (Gaussian, Cauchy, and Laplace), and also under different levels of noise for each type of measurement (range, angle, and Doppler shift). We show under these circumstances that, for the same number of measurements, our method is able to attain the same accuracy as the trilateration.

Under Gaussian and Laplace noise, as the number of radars (therefore, measurements) increases, we see that the estimation error of both position and velocity decreases. Thus, we conclude that our approach can be seen as a generalization of the trilateration method for an arbitrary number of measurements, being able to leverage novel technology to improve accuracy of estimation.

%

\appendices


\subsection{Proof of Proposition 2}
\proptwo*

\paragraph*{Proof.} Remembering the composition of matrix $A_i$ in~\eqref{eq:y_i-quadratic-objective-aux-variables}, we have

\begin{equation}
    \label{eq:matrix-A-reformulated}
    \begin{aligned}
        A_i + \lambda_i^{*} I &= \alpha_i^2 I + \omega_i^2 \beta_i^2 v v^T + \lambda_i^{*} I \\
        &= \eta_i I + \gamma_i v v^T .
    \end{aligned}
\end{equation}
where $\eta_i = \alpha_i^2 + \lambda_i^{*}$ and $\gamma_i = \omega_i^2 \beta_i^2$. A matrix $M$ is invertible if and only if $\det(M) \neq 0$. From the matrix determinant lemma~\cite{ding2007eigenvalues}, we can compute the determinant as

\begin{equation}
    \label{eq:matrix-determinant}
    \begin{aligned}
        \det&\left( \eta_i I + \gamma_i v v^T \right) = \\
        &= \left(1 + v^T \left( \eta_i I \right)^{-1} \left(\gamma_i v\right) \right) \det\left( \eta_i I\right) \\
        &= \left(1 + \frac{\gamma_i}{\eta_i} v^T v\right) \eta_i^n \quad \left(\mbox{in our case, $n = 3$}\right). \\
    \end{aligned}
\end{equation}
Therefore

\begin{equation}
    \label{eq:determinant-condition}
    \begin{aligned}
        \det\left( A_i + \lambda_i^{*} I \right) = 0 &\implies \left(1 + \frac{\gamma_i}{\eta_i} v^T v\right) \eta_i^n = 0 \\
        &\stackrel{\eta_i > 0}{\implies} \|v\|^2 = -\frac{\eta_i}{\gamma_i}.
    \end{aligned}
\end{equation}
Since $\alpha_i^2, \beta_i^2, \omega_i^2 > 0$ and, from~\eqref{eq:dual_feasibility-3}, $\lambda_{i}^{*} \geq 0$, then $-\frac{\eta_i}{\gamma_i} < 0$, therefore $\|v\|^2 \neq -\frac{\eta_i}{\gamma_i}$, which implies that $\left( A_i + \lambda_i^{*} I \right)$ is always invertible for $\lambda_{i}^{*} \geq 0$. \qed

\bibliographystyle{IEEEtaes.bst}
\bibliography{IEEE_TAES_regular_template_latex}

%

\begin{IEEEbiography}[{\includegraphics[width=1in,height=1.25in,clip,keepaspectratio]{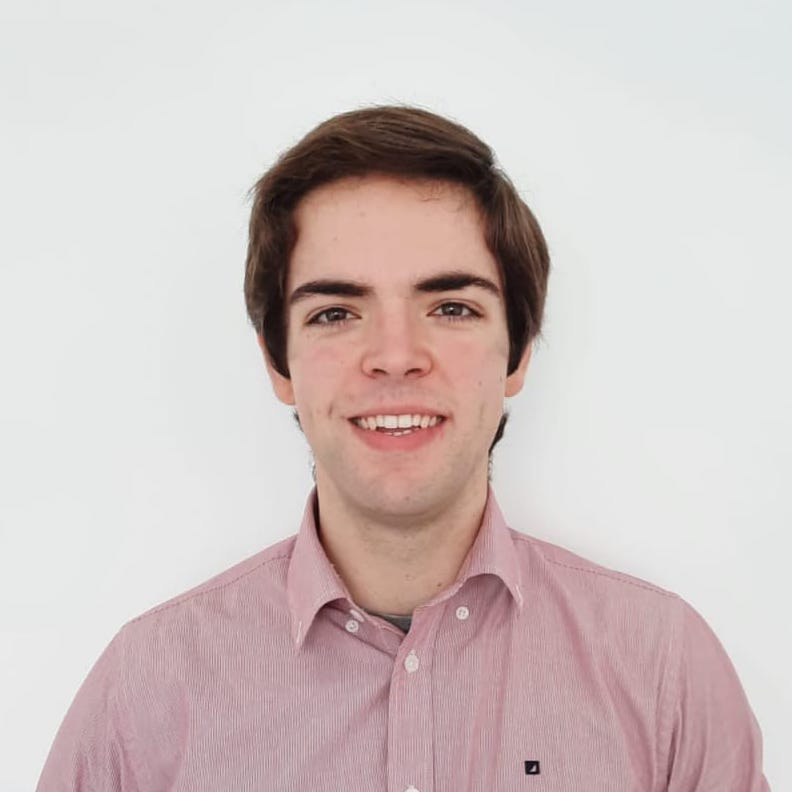}}]{Ricardo Ferreira} received his B.Sc. and M.Sc. degree in Computer Science from NOVA School of Science and Technology. He is currently a Ph.D. student in Computer Science at NOVA School of Science and Technology. His research interests include Machine Learning, Optimization and Probability Theory to develop robust solutions for satellite collision avoidance.
\end{IEEEbiography}%

\begin{IEEEbiography}[{\includegraphics[width=1in,height=1.25in,clip,keepaspectratio]{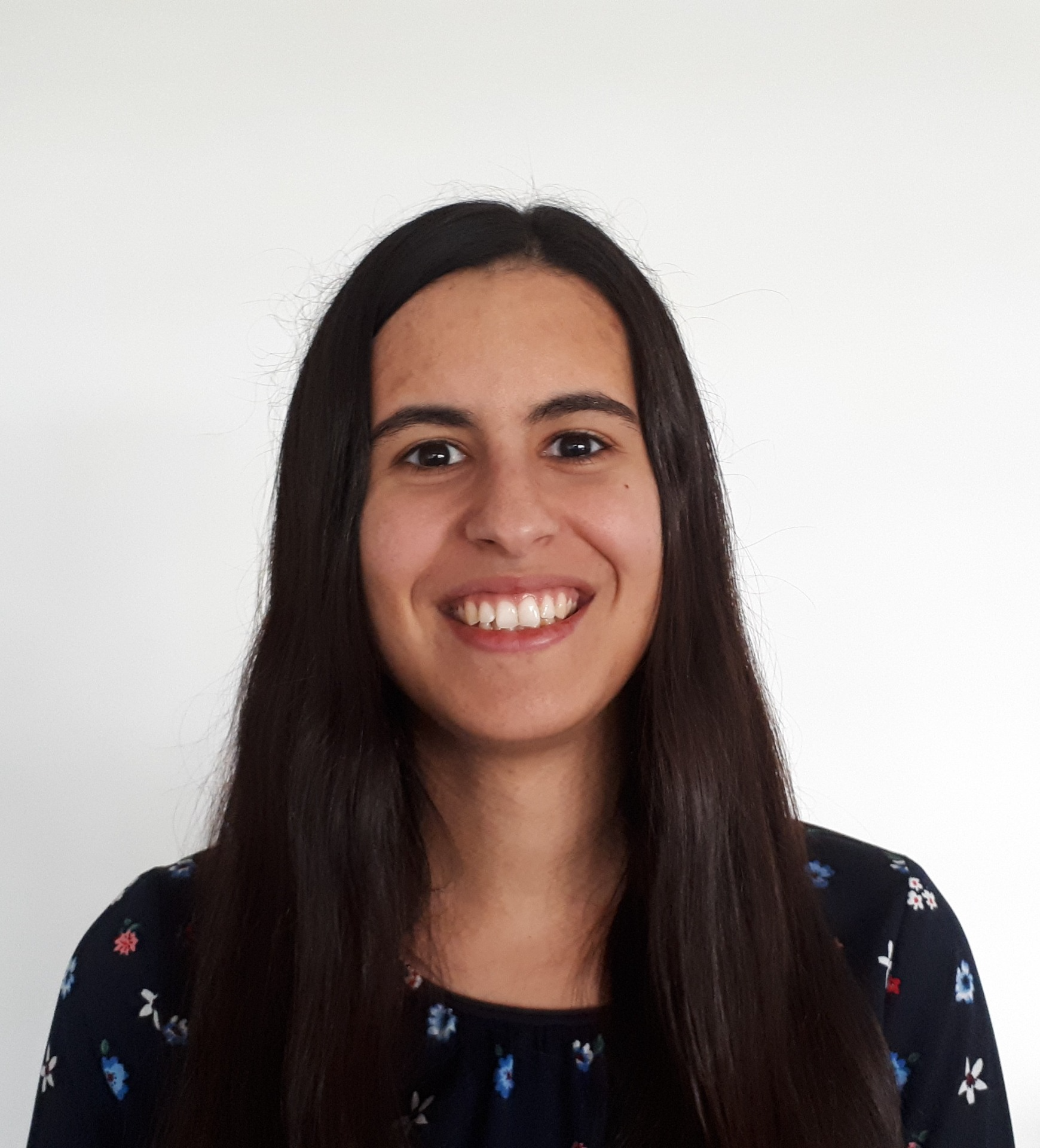}}]{Filipa Valdeira}
received her B.Sc. and M.Sc. in Aerospace Engineering from Instituto Superior Técnico in 2018 and a Ph.D. in Mathematical Sciences from the University of Milan in 2022. She is currently a Post-doctoral Researcher at NOVA School of Science and Technology with NOVA LINCS. Her research interests include Machine Learning, Optimization, 3D Shape Modelling and Gaussian Processes.
\end{IEEEbiography}

\begin{IEEEbiography}[{\includegraphics[width=1in,height=1.25in,clip,keepaspectratio]{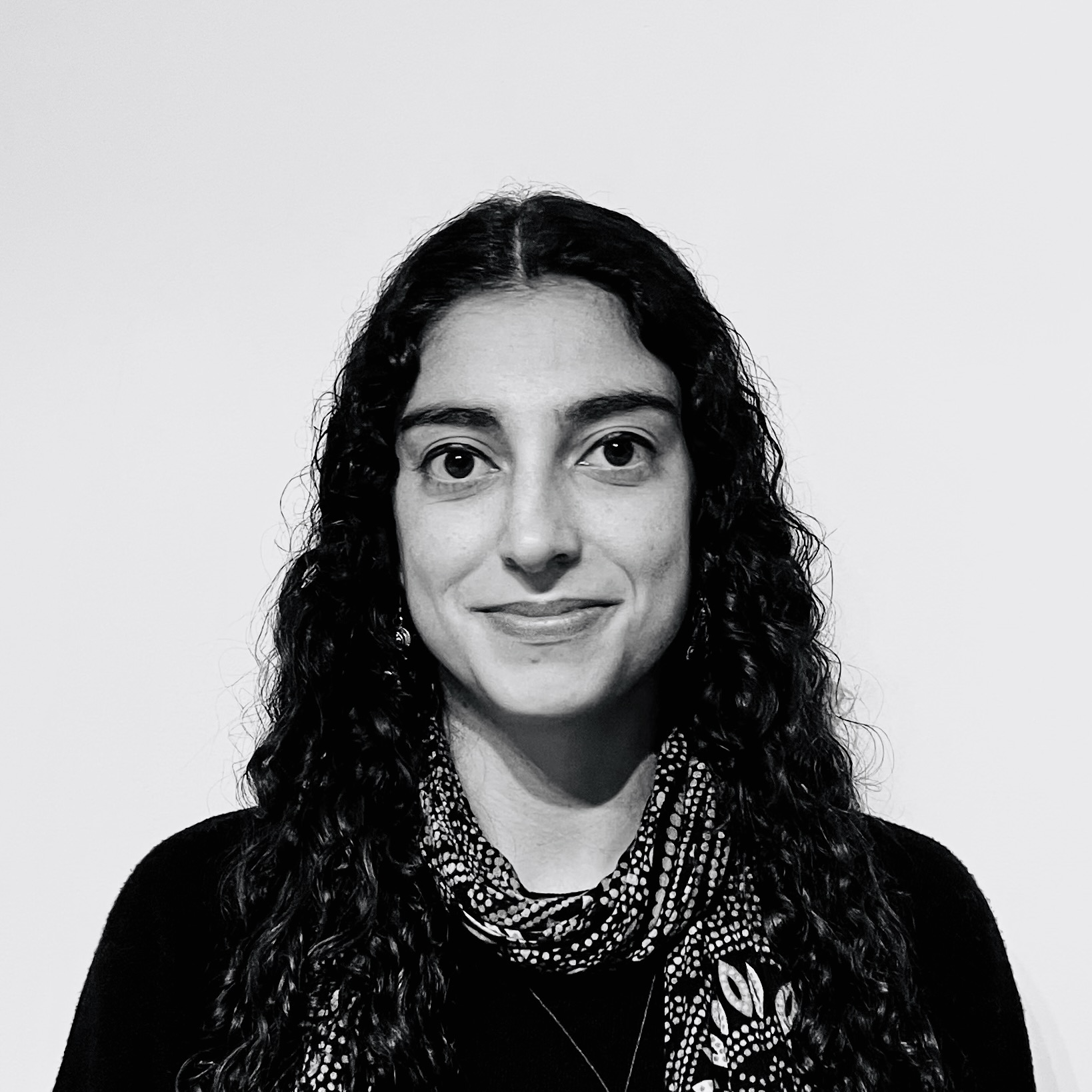}}]{Marta Guimarães} received her B.Sc. and M.Sc. degrees in Aerospace Engineering from Instituto Superior Técnico. She is currently a Ph.D. student in Computer Science at NOVA School of Science and Technology and works as an AI Researcher at Neuraspace, developing Machine Learning solutions for satellite collision avoidance and space debris mitigation.
Her research interests include Machine Learning and Deep Learning, with a focus on Time Series Forecasting.
\end{IEEEbiography}

\begin{IEEEbiography}[{\includegraphics[width=1in,height=1.25in,clip,keepaspectratio]{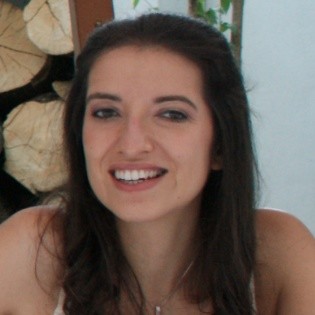}}]{Cláudia Soares} received a Diploma in modern languages and literature from Nova University of Lisbon, Portugal, as well as a B.Sc., M.Sc., and Ph.D. in Engineering from Instituto Superior Tecnico, Portugal. She is currently an Assistant Professor at NOVA School of Science and Technology, Portugal. Her research focuses on using optimization, physics, and probability theory to develop trustworthy machine learning for various applications, including space.
\end{IEEEbiography}

\end{document}